\documentstyle[a4,12pt]{article}

\def\bbbR{{\vbox{\hbox to 8.9pt{I\hskip-2.1pt R\hfil}}}}

\begin{document}

\title{A Brief Survey on Fibrewise General Topology}
\author{Giorgio NORDO}
\date{}
\maketitle

\begin{abstract}
We present some recent results in
Fibrewise General Topology with special regard to the theory of
Tychonoff compactifications of mappings.
Several open problems are also proposed.
\end{abstract}

\parindent=10pt \parskip=3pt

\paragraph{1. Introduction}

\indent

Mapping are more general object than of topological spaces.
In fact, it is evident that any space can be trivially identified with
the continuous mapping of that space to a single-point space.

Since 50's, this simple fact suggested the idea to consider
properties for mappings instead of the traditional ones for spaces
in order to obtain more general statements.

First steps in this direction were moved by Whyburn
\cite{whyb1, whyb2} and Dickman \cite{dickman},
but only in 1975, Ul'janov \cite{uljanov}
introduced the notion of {\em Hausdorff
mapping} (formerly called {\em separable}) to study the Hausdorff
compactifications of countable character.

Later, Pasynkov \cite{pasy} generalized and studied in a systematic way to
the continuous mappings various other notions and properties concerning
spaces like the separation axioms $T_0$, $T_1$, $T_2$, $T_{3\frac 12}$, the
{\em regularity}, the {\em complete regularity}, the {\em normality}, the
{\em compactness} and the {\em local compactness}. A considerable part of
these new definitions and constructions is based on the notion of {\em %
partial topological product} (briefly PTP) introduced and studied by
Pasynkov in \cite{ptp_ann}. The main properties of PTP's, included an
analogous for mappings of the Embedding Lemma, are proved in detail in \cite
{ptp}.

Some weaker separation axioms for continuous mappings such as {\em %
semiregularity} and {\em almost regularity} were introduced and studied in
\cite{p_funz}. The problem of their productivity was investigated in \cite
{prod_pf}.

Let us note that Pasynkov's papers (\cite{pasy}, in particular)
have inspirated James to give a slightly different approach to the
same topic (see \cite{james}).

The generalization to mappings of notions originally defined for spaces
belongs to the more general branch of the {\em Fibrewise General Topology}
(sometimes called {\em General Topology of mappings}) and, from a categorial
point of view, it means to pass from the study of the property of the
category ${\bf Top}$ to those of the category ${\bf Top}_Y$ whose objects
are the continuous mappings into some fixed space $Y$ and whose morphisms
are the continuous functions commutating the triangular diagram of two
objects.

Because a property ${\cal P}_Y$ of ${\bf Top}_Y$ can be considered as a
generalization of some corresponding property ${\cal P}$ of ${\bf Top}$ it
must coincide with ${\cal P}$ when $Y$ is a single-point space, that is when
every object $f:X\to Y$ of ${\bf Top}_Y$ can be identified with the space $X$%
, i.e. with an object of ${\bf Top}$.

In rare case, the analogous ${\cal P}_Y$ in ${\bf Top}_Y$ of some property $%
{\cal P}$ of ${\bf Top}$ is quite evident. For example, to ${\cal P}=$
compactness corresponds ${\cal P}_Y=$ perfectness. In other cases (e.g. for
the separation axioms $T_0$ and $T_1$) the property ${\cal P}_Y$ can be
obtained by requiring that the corresponding property ${\cal P}$ holds on
every fibre of the mappings, but for many properties for mappings, it is
necessary to give new definitions which are more complex than the
corresponding ones for the spaces.

The notion of compactification (i.e. perfect extension) of a continuous
mapping was given first by Whyburn in 1953 \cite{whyb1,whyb2}.

It is worth mentioning that in \cite{perfect} it is presented a
filter based method which allows us to build perfect extension of
every function (not necessarily continuous) between two arbitrary
topological spaces.

However, the first general definition of compactification for mappings
analogous to the well-known notion for spaces, was given by Pasynkov in \cite
{pasy}. In fact, in that paper, using techniques based on PTP's, a method is
described to obtain Tychonoff (i.e. completely regular, $T_0$)
compactifications of Tychonoff mappings between arbitrary spaces, and it is
proved that the poset $TK(f)$ of all the Tychonoff compactifications of a
Tychonoff mapping $f:X\to Y$ admits a maximal compactification $\beta
f:\beta _fX\to Y$ which is the exact analogous, in ${\bf Top}_Y$, of the
Stone--\v {C}ech compactification of a Tychonoff space.

Let us note that K\"{u}nzi and Pasynkov \cite{kunzipas} have completely
described the set $TK(f)$ of all the Tychonoff compactifications of a
Tychonoff mapping $f:X\to Y$ by means of presheaves of the rings $%
C^{*}\left( f^{-1}(U)\right) $ with $U$ open set of $Y$.

Recently, Bludova and Nordo \cite{semilattice} have shown that if a
mapping $f:X\to Y$ is {\em Hausdorff compactifiable} (i.e. it has some Hausdorff
compactification) then there exists the greatest (called "maximal")
compactification $\chi f:\chi _fX\to Y$ in the set $HK(f)$ of all the
Hausdorff compactifications of $f$.

An extension to mappings of the notion of {\em H-closedness}
(see \cite{pw} or \cite{hcerrado} for a complete survey)
was given in \cite{hclosed} by Cammaroto, Fedorchuk and Porter, while some other
generalizations to the mappings of the concepts of realcompactness and
Dieudonn\'e completeness were introduced and studied in \cite{rcompact},
\cite{dieudonne}, \cite{completeness} and \cite{completions}.

The notion of {\em perfect compactification} given by Skljarenko in
\cite{skl} was recently generalized to the mappings in \cite{perfun}
and, in \cite{gandia}, it was proved
that both the maximal realcompactification
$vf$ and the Dieudonn\'e completion $\mu f$ of a Tychonoff mapping $f$
are perfect extensions of $f$.

\paragraph{2. Extension to mappings of notions for spaces}

\indent

Throughout all the paper, the word ''space'' will mean ''topological space''
on which, no separation axiom is assumed and all the mappings will be
supposed continuous unless otherwise specified. If $X$ is a space, $\tau (X)$
will denote the family of all open sets of $X$.

For terms and undefined concepts we refer to \cite{eng}. \vskip 4pt

For any fixed space $Y$, we consider the category ${\bf Top}_Y$ where
\[
Ob({\bf Top}_Y)=\{f\in C(X,Y)\,:\quad X\in Ob({\bf Top})\}
\]
is the class of the {\em objects} and, for every pair $f:X\to Y$, $g:Z\to Y$
of objects,
\[
M(f,g)=\{\lambda \in C(X,Z)\,:\quad g\circ \lambda =f\}
\]
is the class of the {\em morphisms} from $f$ to $g$, whose generic
representant is denoted for short by $\lambda :f\to g$.

A morphism $\lambda : f \to g$ from $f : X \to Y$ to $g : Z \to Y$ will be
called {\em surjective} (resp. {\em closed}, {\em dense}) if $\lambda (X) =
Z $ (resp. $\lambda (X)$ is closed in $Z$, $\lambda (X) $ is dense in $Z$).

If $\lambda : f \to g$ is a surjective morphism, we will say that $g$ is the
{\em image} of $f$ ({\em by the morphism $\lambda$}) and we will write that $%
g = \lambda (f)$.

Moreover, we say that a morphism $\lambda :f\to g$ from $f:X\to Y$ to $%
g:Z\to Y$ is an {\em embedding} (resp. a {\em homeomorphism}) if so is the
function $\lambda :X\to Z$.

A mapping $g:Z \to Y$ is said an {\em extension} of $f:X \to Y$ if there
exists some dense embedding $\lambda : f \to g$ (as usual, we shall identify
$X$ and $f$ by $\lambda (X)$ and $g|_{\lambda (X)}$ respectively).

A morphism $\lambda :g\to h$ between two extensions $g:Z\to Y$ and $h:W\to Y$
of a mapping $f:X\to Y$ will be called {\em canonical} if $\lambda
|_{X}=id_X$. \vskip 6pt

Let us introduce some notions and basic facts about partial products \cite
{ptp_ann,ptp}.

Given two spaces $Y$, $Z$ and an open set $O$ of $Y$, we consider the set $%
P=(Y\backslash O) \cup (O \times Z)$ and the map $p : P \to Y$ defined by $%
p|_{Y\backslash O}=id_{Y\backslash O}$ and $p|_{O\times Z}=pr_O$ where $%
pr_O: O \times Z \to O$ denotes the projection of $O\times Z$ onto $O$.

We will call {\em elementary partial topological product} (briefly EPTP)
{\em of $Y$ and fibre $Z$ relatively to the open set $O$} and we will denote
it by $P(Y,Z,O)$, the space generated on $P$ by the basis ${\cal B}%
(Y,Z,O)=p^{-1}(\tau (Y))\cup \tau (O\times Z)$. \newline
The mapping $p:P\to Y$ above defined will be called the {\em projection} of
the EPTP $P=P(Y,Z,O)$ and it is routine to prove that it is a continuous,
onto, open mapping.

It is evident that the EPTP $P(Y,Z,\emptyset )$ is simply $Y$ and that $%
P(Y,Z,Y)$ coincides with the usual product space $Y\times Z$.
\vskip 4pt

Now, let $Y$ be a topological space, $\{Z_\alpha \}_{\alpha \in \Lambda }$
be a family of spaces and $\{O_\alpha \}_{\alpha \in \Lambda }$ be a family
of open sets of $Y$. For every $\alpha \in \Lambda $, let $P_\alpha
=P(Y,Z_\alpha ,O_\alpha )$ be the EPTP of base $Y$ and fibre $Z_\alpha $
relatively to $O_\alpha $ and $p_\alpha :P_\alpha \to Y$ be its projection.
\newline
We will call {\em partial topological product} (PTP for short) {\em of base $Y$
and fibres $\{Z_\alpha \}_{\alpha \in \Lambda }$ relatively to the open sets
$\{O_\alpha \}_{\alpha \in \Lambda }$ } the {\em fan product} of the spaces $%
\{P_\alpha \}_{\alpha \in \Lambda }$ relatively to the mappings $\{p_\alpha
\}_{\alpha \in \Lambda }$, i.e. the subspace
\[
P=\left\{ t=\left\langle t_\alpha \right\rangle _{\alpha \in \Lambda }\in
\prod_{\alpha \in \Lambda }P_\alpha \,:\quad p_\alpha (t_\alpha )=p_\beta
(t_\beta )\quad \forall \alpha ,\beta \in \Lambda \right\}
\]
and we will denote it by $P(Y,\{Z_\alpha \},\{O_\alpha \};\alpha \in \Lambda
)$.

For every $\alpha \in \Lambda$, the restriction $\pi_\alpha = pr_\alpha |_P
: P \to P_\alpha$ of the $\alpha$-th canonical projection $pr_\alpha$ will
be called the $\alpha$-th {\em short projection}, while the {\em fibrewise
product of the mappings} $\{ p_\alpha \}_{\alpha \in \Lambda}$, i.e. the
continuous mapping $p : P \to Y$ defined by $p_\alpha \circ \pi_\alpha = p$
for any $\alpha \in \Lambda$ will be said the {\em long projection} of the
PTP $P(Y, \{ Z_\alpha \} , \{ O_\alpha \} ;\alpha \in \Lambda )$.

In case $O_\alpha = Y$ for every $\alpha \in \Lambda$, the
PTP $P(Y, \{ Z_\alpha \} , \{ O_\alpha \} ;\alpha \in \Lambda )$ coincides
(up to homeomorphisms) with the usual product $Y \times \prod_{\alpha \in
\Lambda} Z_\alpha$, if $|O_\alpha |=1$ for every $\alpha \in \Lambda$, the
PTP $P(Y, \{ Z_\alpha \} , \{ O_\alpha \} ;\alpha \in \Lambda )$ coincides
with the usual Tychonoff product $\prod_{\alpha \in \Lambda} Z_\alpha$ of
its fibres, while if $O_\alpha = \emptyset$ for any $\alpha \in \Lambda$,
the PTP $P(Y, \{ Z_\alpha \} , \{ O_\alpha \} ;\alpha \in \Lambda ) $ is
simply (homeomorphic to) the space $Y$. \vskip 3pt

\noindent {\bf Definitions.} A mapping $f:X\to Y$ is said to be $T_0$ \cite
{pasy} if for every $x,x^{\prime }\in X$ such that $x\ne x^{\prime }$ and $%
f(x)=f(x^{\prime })$ there exists a neighborhood of $x$ in $X$ which does
not contain $x^{\prime }$ or a neighborhood of $x^{\prime }$ in $X$ not
containing $x$.

A mapping $f:X\to Y$ is said to be {\em Hausdorff} (or $T_2$) \cite
{uljanov,pasy} if for every $x,x^{\prime }\in X$ such that $x\ne x^{\prime }$
and $f(x)=f(x^{\prime })$ there are two disjoint neighborhoods of $x$ and $%
x^{\prime }$ in $X$.

We will say that $f: X \to Y$ is {\em compact} if it is perfect (i.e. it is
closed and every its fibre is compact).

A mapping $f:X\to Y$ is said to be {\em completely regular} \cite{pasy} if
for every closed set $F$ of $X$ and $x\in X\backslash F$ there exists a
neighborhood $O$ of $f(x)$ in $Y$ and a continuous function $\varphi
:f^{-1}(O)\to [0,1]$ such that $\varphi (x)=0$ and $\varphi (F\cap
f^{-1}(O))\subseteq \{1\}$.

A completely regular, $T_0$ mapping is called {\em Tychonoff} (or $%
T_{3\frac 12}$) \cite{pasy}. \vskip 4pt

\noindent {\bf Remark.} It is easy to verify that all the previous
properties in ${\bf Top}_Y$ coincide with the corresponding ones in ${\bf Top%
}$ provided $|Y|=1$ and that every continuous mapping $f:X\to Y$ has such a
property $iff$ both the spaces $X$ and $Y$ have the corresponding properties
(in particular, they are ${\cal P}-$functions in the sense of \cite{p_funz}%
). \vskip 5pt

\noindent {\bf Definition.} A restriction $f|_{X^{\prime}}: X^{\prime}\to Y$
to $X^{\prime}\subseteq X$ of a mapping $f : X \to Y$ is said a {\em closed
restriction} of $f$, if $X^{\prime}$ is a closed subset of $X$. \vskip 4pt

Obviously (see for example \cite{pw}), every closed restriction of a compact
mapping is compact too.
\vskip 3pt

Most well-known statements which hold in the category ${\bf Top}$ have
correspondent ones (and hence generalizations) in ${\bf Top}_Y$.
The following properties are essentially given in \cite{pasy}
(detailed proofs can be found in \cite{tesi}).

\proclaim PROPOSITION 2.1.
Every image $\lambda (f)$ of a compact mapping $f:X\to Y$ is compact too.

\proclaim PROPOSITION 2.2.
Every closed restriction $f|_{X'}$ of a compact mapping $f:X\to Y$
is compact too.

\proclaim PROPOSITION 2.3.
Every compact restriction $f|_{X'}$ of a Hausdorff mapping $f:X \to Y$
is a closed restriction of $f$.

\proclaim PROPOSITION 2.4. Let $\lambda $ and $\mu $ be morphisms from a
mapping $f:X\to Y$ to a Hausdorff mapping $g:Z\to Y$ and $D$ be a dense
subset of $X$. Then if $\lambda |_D =\mu |_D$, the morphisms $\lambda $ and $%
\mu $ coincide.

\proclaim PROPOSITION 2.5. Every morphism $\lambda : f \to g$ from a
compact mapping $f : X \to Y$ to a Hausdorff mapping $g : Z \to Y$ is
perfect.

\paragraph{3. Compactification of mappings}

\indent

Let $f:X \to Y$ be a mapping.
We say that a mapping $c:X^c\to Y$ is a
{\em compactification} of $f$ (in ${\bf Top}_Y$) if it is a compact
($=$ perfect) extension of $f$.\vskip 2pt

This approach to the notion of the compactification
of a mapping was proposed by Whyburn \cite{whyb2}, but in the most general
situation, this notion was studied first by Pasynkov in \cite{pasy}.
\vskip 4pt

\noindent
{\bf Remark.} A different variant of
compactifications of mappings was examined by Uljanov \cite{uljanov}.
But, it is a common opinion that Uljanov's definition is not natural for
non-surjective mappings because, in that case, a compact mapping is
not its own compactification.
\vskip 5pt

\noindent {\bf Definitions.} Let $c:X^c\to Y$ and $d:X^d\to Y$ be two
compactifications of a mapping $f:X\to Y$ (in ${\bf Top}_Y$). We say that: %
\parskip=-5pt
\begin{itemize}
\item  {\em $c$ is projectively larger than $d$} (relatively to $f
$) and we write that $c\ge _fd$ (or $c\ge d$, for short) if there exists
some canonical morphism $\lambda :c\to d$.
\item  {\em $c$ is equivalent to $d$} (relatively to $f$) and we
write that $c\equiv _fd$ (shortly, $c\equiv d\,$) if there exists a
canonical homeomorphism $\lambda :c\to d$.
\end{itemize}
\parskip=3pt

The following useful result is given in \cite{semilattice}.

\proclaim PROPOSITION 3.1. Let $c:X^c\to Y$ and $d:X^d\to Y$ be two
Hausdorff compactifications of a mapping $f:X\to Y$. Then $c\equiv _fd$ if
and only if $c\ge_f d$ and $d\ge_f c$.

\noindent {\bf Definition.} A Hausdorff mapping $f:X\to Y$ will be called
{\em Hausdorff compactifiable} if it has some Hausdorff compactification (in ${\bf Top}%
_Y$). \vskip 3pt

In \cite{semilattice}, it is noted that the class of all Hausdorff
compactifications of any Hausdorff compactifiable mapping $f:X\to Y$ forms a set
modulo the equivalence $\equiv _f$. \vskip 2pt

\noindent {\bf Definition.} If $f : X \to Y$ is a Hausdorff compactifiable mapping, $%
HK(f)$ will denote the set of all equivalence classes of Hausdorff
compactifications of $f$. \vskip 3pt

So, by 3.1, it follows that $(HK(f),\ge )$ is a poset and, for any pair of
Hausdorff compactifications $c,d\in HK(f)$ we can write $c=d$ instead of $%
c\equiv _fd$, that is we do not distinguish between equivalent Hausdorff
compactifications. \vskip 3pt

In \cite{pasy}, Pasynkov erroneously indicated that it is proved
in \cite{uljanov} that every Hausdorff compactifiable mapping $f:X\to Y$
has a maximal one.
\\
This fact was also used in several following papers like
\cite{dieudonne}, \cite{hclosed}, \cite{rcompact},
\cite{kunzipas}, \cite{completeness}, \cite{mazroa},
\cite{completions}, etc. but it is not correct because
the Ul'janov's definition is different from the currently used one
(given by Pasynkov in \cite{pasy}) that does not include the surjectivity.
\\
Anyway the existence of the maximal Hausdorff compactification
was actually proved by Bludova and the author as direct consequence
of the following more general result.
\vskip 6pt

\noindent {\bf THEOREM 3.2.} \cite{semilattice}
{\sl For any Hausdorff compactifiable mapping $f:X\to Y$, $%
(HK(f),\ge )$ is a complete upper semilattice}
\vskip 8pt

The projective maximum of $(HK(f), \ge )$, i.e. the {\em maximal Hausdorff
compactification} of $f$, will be denoted by $\chi f : \chi _f X \to Y$. %
\vskip 2pt

From this and by 2.4 and 2.5, it follows -- in particular -- that for
any Hausdorff compactification $bf:X^b\to Y$ of a Hausdorff
compactifiable mapping $f:X\to Y$
there exists a unique perfect canonical morphism $\lambda _b:\chi
f\to bf$.

In \cite{pasy}, Pasynkov proved that every Tychonoff mapping has a
Tychonoff compactification.

Since it is easy to show that a Tychonoff mapping is Hausdorff,
Proposition 3.1 allow us to give the following: \vskip 2pt

\noindent {\bf Definition.} For any Tychonoff mapping $f:X\to Y$, we will
denote by $TK(f)$ the set of all Tychonoff compactifications of $f$ up to
the equivalence $\equiv _f$.
\vskip 3pt

For a mapping $f:X \to Y$, let us denote by $C^*(f)$ the family of all
the {\em partial mappings} on $f$, i.e. of all the continuous bounded
real-valued mappings $\varphi : f^{-1} (O_\varphi) \to I_\varphi$ defined
from the inverse image by $f$ of an open set $O_\varphi$ of $Y$
to a compact subset $I_\varphi$ of the real line $\bbbR$.
\vskip 3pt

\noindent
{\bf Definitions.}
A subfamily ${\cal C }= \{ \varphi :f^{-1}(O_\varphi)
\to I_\varphi \}$ of $C^*(f)$
is said to be:
\parskip=-6pt
\begin{itemize}
\item {\em separating the points of $f$} if for every $x,x' \in X$
such that $f(x)=f(x')$ there exists some $\varphi \in {\cal C}$
such that $x,x' \in f^{-1} (O_\varphi )$ and $\varphi (x) \ne \varphi (x')$.
\vspace{-3mm}
\item {\em separating the points from the closed sets of $f$}
if for any closed set $F$ of $X$ and every $x \in X \setminus F$
there exists some $\varphi \in {\cal C}$
such that $x \in f^{-1} (O_\varphi )$ and
$\varphi (x) \notin cl_{I_\varphi } \left( F
\cap f^{-1} (O_\varphi ) \right)$.
\end{itemize}
\vskip 6pt

It is shown in \cite{pasy} (see \cite{tesi} for a more detailed proof)
that every Tychonoff compactification $bf:X^b \to Y$ of a Tychonoff mapping
$f:X\to Y$ is uniquely determinated by a subfamily
${\cal C} = \{ \varphi: f^{-1}(O_\varphi ) \to
I_\varphi  \, : \, \, O_\varphi \in \tau (Y) \}$
of $C^*(f)$ separating the points and the points from the closed sets of $f$
and that $bf$ coincides with a particular restriction of the long projection
$p_{\cal C} : P_{\cal C} \to Y$ of the PTP
$P_{\cal C} = P( Y , \{ O_\varphi \} ,
\{ I_\varphi \} ; \varphi \in {\cal C } )$.
\vskip 12pt

Thus, the notion of PTP plays in the category ${\bf Top}_Y$ the same role
that the notion of product space has in the category ${\bf Top}$ and, as
matter of fact, they coincide when $|Y|=1$.
\vskip 10pt

In \cite{pasy}, it is also proved that for any Tychonoff mapping $f:X\to Y$
there exists, in $(TK(f),\ge )$, a maximal Tychonoff compactification $\beta
f:\beta _fX\to Y$ that is determinated by all the whole family $C^*(f)$
and characterized by some extension properties very similar to that of
the Stone-\v{C}ech compactification.
\\
We have, in fact, the following:
\vskip 12pt

\proclaim THEOREM 3.3. For any Tychonoff compactification $bf:X^b\to Y$ of a
Tychonoff mapping $f:X\to Y$, the following conditions are equivalent: %
\parskip=-6pt
\begin{enumerate}
\item[(1)]  $bf=\beta f$ ;
\vspace{-3mm}
\item[(2)]  for every $U\in \tau (Y)$ and $\varphi \in C^{*}\left(
f^{-1}(U)\right) $ there exists a unique extension $\widetilde{\varphi }\in
C^{*}\left( (bf)^{-1}(U)\right) $ ;
\vspace{-3mm}
\item[(3)]  for every compact Tychonoff mapping $k:Z\to Y$ and every
morphism $\lambda :f\to k$ there exists a morphism $\widetilde{\lambda }%
:bf\to k$ which extends $\lambda $.
\end{enumerate}

\parskip=3pt

Moreover, Theorem 3.2. allow us to obtain as immediate consequence the
following:
\vskip 4pt

\noindent {\bf THEOREM 3.4.} \cite{semilattice} {\sl The poset $(TK(f),\ge )$
of all Tychonoff compactifications of a Tychonoff mapping $f:X\to Y$ is a
complete upper semilattice whose projective maximum is $\beta f$.}
\vskip 6pt

\noindent {\bf PROPOSITION 3.5.} \cite{pasy} {\sl For any Tychonoff
compactification $bf:X^b\to Y$ of a Tychonoff mapping $f:X\to Y$ there
exists a unique (perfect) canonical morphism $\mu _b:\beta f\to bf$ such
that $\mu _b(\beta _fX\backslash X)=X^b\backslash X$.}
\vskip 6pt

Let us observe that if $|Y|=1$, $X$ is a Tychonoff space, the domain $\beta
_fX$ of $\beta f$ coincides with the Stone \v {C}ech compactification $\beta
X$ of $X$, the domain $X^b$ of $bf$ is a generic compactification of $X$ and
$\lambda _b:\beta _fX\to X^b$ becomes the usual quotient map (see for
example \cite{chandler}). \vskip 5mm

In general, for a Tychonoff mapping $f:X\to Y$, we have
\[
TK(f)
\begin{array}{c}
\\
\subset  \\
\ne
\end{array}
HK(f)
\]
that is, unlike the corresponding case for spaces, there exist Hausdorff
compactification which are not Tychonoff or, equivalently, there are compact
Hausdorff mapping which are not Tychonoff. In fact, it was proved in \cite
{chaber} (see also \cite{isbell}) that it is possible to build a perfect ($%
\equiv $ compact) mapping defined on a regular $T_0$ but non Tychonoff
space onto a Tychonoff space (that is the property $T_{3\frac 12}$ is not an
inverse invariant by perfect mappings) and since it is proved in \cite{pasy}
that if a mapping and its range are both completely regular, its domain is
too, it follows directly that such a mapping can not be completely regular
and hence Tychonoff.

This is the reason why it is necessary to study the classes of Tychonoff
and Hausdorff compactifiable mappings separately.

\paragraph{4. Open problems.}

\indent

It seems that the following questions might be interesting.
Some of these problems are published for the first time.
\vskip 5pt

\noindent
{\bf Problem 1.}
It is well-known that the poset $K(X)$
of all Hausdorff compactification of a Tychonoff
space $X$ can be completely characterized in terms
of the families of $C^*(X)$ that separete points and points
from closed sets of $X$ (see, for example, \cite{chandler}).
\\
{\em Is it possible to obtain such a similar characterization for
the set $HK(f)$ (the set $TK(f)$) of all Hausdorff (Tychonoff)
compactification of a Hausdorff compactifiable (Tychonoff) mapping $f$
in terms of the separating families of $C^*(f) \,$?}
\vskip 4pt

\noindent
{\bf Problem 2.}
Magill has proved in \cite{magill} (see also \cite{chandler})
that the posets $K(X)$ and $K(Y)$
of all Hausdorff compactifications of two locally compact spaces $X$ and $Y$
are isomorphic if and only if their Stone-\v{C}ech remainders
$\beta X \setminus X$ and $\beta Y \setminus Y$ are homeomorphic.
\\
{\em Is it possible to find a definition of locally compact mapping
that allow us to obtain a Magill-type theorem for mappings$\,$?}
\vskip 4pt

\noindent
{\bf Problem 3.}
{\em Is it possible to obtain analogous in ${\bf Top}_Y$
of properties like the countably compactness,
the paracompactness and the pseudocompactness$\,$?}
\vskip 4pt

\noindent
{\bf Problem 4.}
{\em Is there a consistent definition of metrizable mapping
which extends the corresponding notion for spaces
and allow us to obtain general metrization theorems$\,$?}
\vskip 4pt

\noindent
{\bf Problem 5.}
{\em Is it possible to extend to mappings
other kind of Tychonoff extension properties
like the $m$-boundedness {\rm(see \cite{pw})} ?}

\vspace{-1mm}

{\noindent {\it Key words and phrases:}} partial topological product,
$T_0$ mapping, Hausdorff mapping, completely regular
mapping, Tychonoff mapping, compact mapping, compactification of a
mapping, realcompact mapping, Dieudonn\'{e} complete mapping.
\vskip 6pt

{\noindent {\it AMS Subject Classification:}} Primary 54C05, 54C10, 54C20,
54C25; Secondary: 54D15, 54D35, 54D60.
\vskip 6pt

\parskip=0pt {\noindent {\sc Giorgio NORDO}\newline
MIFT - Dipartimento di Scienze Matematiche e Informatiche, scienze Fisiche e scienze della Terra,
Messina University, Messina, Italy
\vskip 3pt \noindent E-mail: {\tt giorgio.nordo@unime.it}


\parskip=0pt
\begin{thebibliography}{CGNP}
\vspace{-2mm}
\bibitem[BN]{semilattice}  BLUDOVA~I.V., NORDO~G., {\it On the posets of all
the Hausdorff and all the Tychonoff compactifications of continuous mappings}%
, Q \& A in General Topology, Vol. {\bf 17} (1999), 47-55.

\bibitem[BuP]{dieudonne}  BUZULINA~T.I., PASYNKOV~B.A., {\it On
Dieudonn\'{e} complete mappings}, in: Geometry of immersed manifolds, Izdat.
''Prometei'' MGPI, Moscow (1989), 95-98 (in Russian).

\bibitem[CFP]{hclosed}  CAMMAROTO~F., V.V.~FEDORCHUK, J.R.~PORTER, {\it On
H-closed functions}, Comment. Math. Univ. Carolinae {\bf 39},3 (1998),
563-572.

\bibitem[CGNP]{hcerrado} CAMMAROTO~F, GUTIERREZ~J., NORDO~G., de~PRADA, M.A.,
{\it Introduccion a los espacios H-cerrados -- Principales
contribuciones a las formas debiles de compacidad -- Problemas
abiertos}, Mathematic\ae$\,$ Not\ae$\,$ {\bf 38} (1995-96), 47-77 (in Spanish).

\bibitem[CN]{p_funz}  CAMMAROTO~F.,~NORDO~G., {\it On Urysohn, almost
regular and semiregular functions}, Filomat n. {\bf 8} (1994), 71-80.

\bibitem[Cb]{chaber}  CHABER~J., {\it Remarks on open-closed mappings},
Fund. Math. {\bf 74} (1972), 197-208.

\bibitem[Ch]{chandler}  CHANDLER~R.E., {\it Hausdorff compactifications},
Marcel Dekker, New York, 1976.

\bibitem[D]{dickman}
R.F. DICKMAN~Jr., {\it On closed extensions of functions},
Proc. Nat. Acad. Sci. U.S.A. {\bf 62} (1969), 326-332.

\bibitem[E]{eng}  ENGELKING~R., {\it General Topology}, Heldermann, Berlin,
1989.

\bibitem[HI]{isbell}  HENRIKSEN~M., ISBELL~J.R., {\it Some properties of
compactifications}, Duke Math. Journal {\bf 25} (1958), 83-106.

\bibitem[IP]{rcompact}  IL'INA~N.I., PASYNKOV~B.A., {\it On R--complete
mappings}, in: Geometry of immersed manifolds, Izdat. ''Prometei'' MGPI,
Moscow (1989), 125-130 (in Russian).

\bibitem[J]{james}
I.M. JAMES, {\it Fibrewise Topology},
Cambridge University Press, Cambridge, 1989.

\bibitem[KP]{kunzipas}  K\"{U}NZI~H.P.A.,~PASYNKOV~B.A., {\it Tychonoff
compactifications and R-completions of mappings and rings of continuous
functions}, Applied Categor. Structures, {\bf 4} (1996), 175-202.

\bibitem[M]{magill}  MAGILL~K.D., {\it The lattice of compactifications
of a locally compact space}, Proc. London Math. Soc. {\bf 18} (1968),
231-244.

\bibitem[MuP]{completeness}  MUSAEV~D.K., PASYNKOV~B.A., {\it On properties
of compactness and completeness of topological spaces and continuous mappings%
}, in: Tashkeit FAN, Acad. of Sci. of Republic Uzbekistan, 1994 (in Russian).

\bibitem[M$_2$]{mazroa}  MAZROA~E.M.R., {\it Perfect compactifications of
continuous mappings}, Vestn. Mosk. Univ. Ser. I (1990), no. {\bf 1}, 23-26.
= Moscow Univ. Math. Bull. {\bf 45} (1990), no. 1, 24-26.

\bibitem[N$_1$]{prod_pf}  NORDO~G., {\it On product of $\,{\cal P}-$functions%
}, Atti Accad. Peloritana dei Pericolanti, Classe I di Scienze MM.FF.NN.
Vol. {\bf LXXII} (1994), 465-478.

\bibitem[N$_2$]{perfect}  NORDO~G., {\it A note on perfectification of
mappings}, Q \& A in General Topology, Vol. {\bf 14} (1996), 107-110.

\bibitem[N$_3$]{perfext}  NORDO~G., {\it A basic approach to the perfect
extensions of spaces}, Comment. Math. Univ. Carolinae {\bf 38},3 (1997),
571-580.

\bibitem[N$_4$]{tesi}  NORDO~G., {\it Compattificazioni perfette di funzioni},
Ph.D. Dissertation, Messina, 1998 (in Italian).

\bibitem[NP]{perfun}  NORDO~G., PASYNKOV~B.A., {\it Perfect
compactifications of functions}, Comment. Math. Univ. Carolinae {\bf 41},3 (2000),
619-529.

\bibitem[No]{norin}  NORIN~V.P., {\it On proximities for mappings}, Vestn.
Mosk. Univ. Math. Mech. Ser. I 1982, no. {\bf 4} (1982), 33-36 = Moscow
Univ. Math. Bull. {\bf 37}, no. 4 (1982), 40-44.

\bibitem[P$_1$]{ptp_ann}  PASYNKOV~B.A., {\it Partial topological products},
Akad. Nauk S.S.S.R., {\bf 154} (1964), 767-770.

\bibitem[P$_2$]{ptp}  PASYNKOV~B.A., {\it Partial topological products},
Trudy Moskov. Matem. Obshchestva {\bf 13} (1965), 136-245 = Trans. Moscow
Math. Soc. {\bf 13} (1965), 153-272.

\bibitem[P$_3$]{pasy}  PASYNKOV~B.A., {\it On extension to mappings of
certain notions and assertions concerning spaces}, in: Mapping and Functors,
Izdat. MGU, Moscow (1984), 72-102 (in Russian).

\bibitem[P$_4$]{completions}  PASYNKOV~B.A., {\it On completions of mappings}%
, in: Geometry of immersed manifolds, Izdat. ''Prometei'' MGPI, Moscow
(1989), 131-136 (in Russian).

\bibitem[PoW]{pw}  PORTER~J.R., WOODS~R.G., {\it Extensions and absolutes of
Hausdorff spaces}, Springer, 1988.

\bibitem[S]{skl}  SKLJARENKO~E.G., {\it On perfect bicompact extensions},
Dokl. Akad. Nauk S.S.S.R. {\bf 137} (1961), 39-41 = Soviet Math. Dokl. {\bf 2%
} (1961), 238-240.

\bibitem[U]{uljanov}  UL'JANOV~V.M., {\it On compactifications satisfying
the first axiom of countability and absolutes}, Math USSR Sbornik, Vol. {\bf %
27}, n.2 (1975), 199-226.

\bibitem[W$_1$]{whyb1}  WHYBURN~G.T., {\it A unified space for mappings},
Trans.~A.M.S. {\bf 74} (1953), 344-350.

\bibitem[W$_2$]{whyb2}  WHYBURN~G.T., {\it Compactification of mappings},
Math.~Ann. {\bf 166} (1966), 168-174.
\end{thebibliography}
\end{document}